\documentstyle[11pt]{article}
\begin{document}
\renewcommand{\thefootnote}{\fnsymbol{footnote}}
\pagestyle{empty} \setcounter{page}{1}



\def\id{{\bf 1}}
\def\np{\vfil\eject}
\def\nl{\hfil\break}

\newfont{\twelvemsb}{msbm10 scaled\magstep1}
\newfont{\eightmsb}{msbm8} \newfont{\sixmsb}{msbm6} \newfam\msbfam
\textfont\msbfam=\twelvemsb \scriptfont\msbfam=\eightmsb
\scriptscriptfont\msbfam=\sixmsb \catcode`\@=11
\def\Bbb{\ifmmode\let\next\Bbb@\else \def\next{\errmessage{Use
      \string\Bbb\space only in math mode}}\fi\next}
\def\Bbb@#1{{\Bbb@@{#1}}} \def\Bbb@@#1{\fam\msbfam#1}
\newfont{\twelvegoth}{eufm10 scaled\magstep1}
\newfont{\tengoth}{eufm10} \newfont{\eightgoth}{eufm8}
\newfont{\sixgoth}{eufm6} \newfam\gothfam
\textfont\gothfam=\twelvegoth \scriptfont\gothfam=\eightgoth
\scriptscriptfont\gothfam=\sixgoth \def\frak{\frak@}
\def\frak@#1{{\fam\gothfam{{#1}}}} \def\frak@@#1{\fam\gothfam#1}
\catcode`@=12

\def\CC{{\Bbb C}}
\def\NN{{\Bbb N}}
\def\QQ{{\Bbb Q}}
\def\RR{{\Bbb R}}
\def\ZZ{{\Bbb Z}}
\def\cA{{\cal A}}          \def\cB{{\cal B}}          \def\cC{{\cal C}}
\def\cD{{\cal D}}          \def\cE{{\cal E}}          \def\cF{{\cal F}}
\def\cG{{\cal G}}          \def\cH{{\cal H}}          \def\cI{{\cal I}}
\def\cJ{{\cal J}}          \def\cK{{\cal K}}          \def\cL{{\cal L}}
\def\cM{{\cal M}}          \def\cN{{\cal N}}          \def\cO{{\cal O}}
\def\cP{{\cal P}}          \def\cQ{{\cal Q}}          \def\cR{{\cal R}}
\def\cS{{\cal S}}          \def\cT{{\cal T}}          \def\cU{{\cal U}}
\def\cV{{\cal V}}          \def\cW{{\cal W}}          \def\cX{{\cal X}}
\def\cY{{\cal Y}}          \def\cZ{{\cal Z}}
\def\qed{\hfill \rule{5pt}{5pt}}
\def\arcsinh{\mathop{\rm arcsinh}\nolimits}
\newtheorem{theorem}{Theorem}
\newtheorem{prop}{Proposition}
\newtheorem{conj}{Conjecture}
\newenvironment{result}{\vspace{.2cm} \em}{\vspace{.2cm}}
\pagestyle{plain}
\begin{center}

{\LARGE {\bf ADDENDUM TO "ON SUPER-JORDANIAN ${\cal U}_{\sf
h}(sl(N|1))$ ALGEBRA"}}
\\[0.2cm]

\smallskip

{\small B. ABDESSELAM$^{a}$, A. CHAKRABARTI$^{b}$, R.
CHAKRABARTI$^{c}$, A. YANALLAH$^{d}$ and M.B ZAHAF$^{e}$}

\smallskip

{\small {\it $^{a,d,e}$Laboratoire de Physique Quantique de la
Mati\`ere et Mod\'elisations Math\'ematiques (LPQ3M), Centre
Universitaire de Mascara, 29000-Mascara, Alg\'erie}}

\smallskip

{\small {\it $^{a}$Laboratoire de Physique Th\'eorique d'Oran,
Universit\'e d'Oran Es-S\'enia, 31100-Oran, Alg\'erie}}

\smallskip

{\small {\it $^{b}$Centre de Physique Th\'eorique, Ecole
Polytechnique, 91128-Palaiseau cedex, France.}}

\smallskip

{\small {\it $^{c}$Department of Theoretical Physics, University
of Madras, Guindy Campus, Madras 600025, India}}

\end{center}

\begin{abstract}

\noindent We give a complete proof of the result (2.10) presented
in our paper published in {\bf J. Phys. A: Math. Gen. 39 (2006)
8307–--8319}.

\end{abstract}

1. Let ${\cal U}_q(sl(2|1))$ ($q$ is an arbitrary complex number)
be the Hopf superalgebra generated by the elements ${\hat h}_i$,
${\hat e}_i$ and ${\hat f}_i$, $i=1,2$, satisfy the relations
\begin{eqnarray}
&&[{\hat h}_i,\;{\hat h}_j]=0,\qquad \qquad [{\hat h}_i,\;{\hat
e}_j]=a_{ij}{\hat e}_j,\qquad \qquad
[{\hat h}_i,\;{\hat f}_j]=-a_{ij}{\hat f}_j,\nonumber\\
&&[{\hat e}_i,\;{\hat f}_j]=\delta_{ij}\frac{q^{{\hat
h}_i}-q^{-{\hat h}_i}}{q-q^{-1}},\qquad
\qquad {\hat e}_2^2={\hat f}_2^2=0,\nonumber\\
&&{\hat e}_1^2{\hat e}_2-\left(q+q^{-1}\right){\hat e}_1{\hat
e}_2{\hat e}_1+ {\hat e}_2{\hat e}_1^2=0,\qquad {\hat f}_1^2{\hat
f}_2-\left(q+q^{-1}\right){\hat f}_1{\hat f}_2{\hat f}_1+ {\hat
f}_2{\hat f}_1^2=0,
\end{eqnarray}
where $[\;,\;]$ is the supercommutator given by
$[a,\;b]=ab-\left(-\right)^{\deg(a)\deg(b)}ba$ and $a_{11}=2$,
$a_{12}=a_{21}=-1$, $a_{22}=0$. All generators are even except for
${\hat e}_2$ and ${\hat f}_2$ which are odd. The coproducts,
counits and antipodes are given by
\begin{eqnarray}
&&\Delta\left({\hat e}_i\right)={\hat e}_i\otimes q^{{\hat
h}_i/2}+q^{-{\hat h}_i/2}\otimes {\hat
e}_i,\qquad\qquad\epsilon({\hat e}_i)=0,\qquad\qquad S({\hat
e}_i)=-q^{{\hat h}_i/2}{\hat e}_iq^{-{\hat h}_i/2},
\nonumber\\
&&\Delta\left({\hat f}_i\right)={\hat f}_i\otimes q^{{\hat
h}_i/2}+q^{-{\hat h}_i/2}\otimes {\hat
f}_i,\qquad\qquad\epsilon({\hat f}_i)=0,\qquad\qquad S({\hat
f}_i)=-q^{{\hat h}_i/2}{\hat f}_iq^{-{\hat h}_i/2},
\nonumber\\
&&\Delta\left({\hat h}_i\right)={\hat h}_i\otimes 1+1\otimes {\hat
h}_i,\qquad\qquad\epsilon({\hat h}_i)=0,\qquad\qquad S({\hat
h}_i)=-{\hat h}_i.
\end{eqnarray}
Let ${\hat e}_3={\hat e}_1{\hat e}_2-q^{-1}{\hat e}_2{\hat e}_1$
and ${\hat f}_3={\hat f}_2{\hat f}_1-q{\hat f}_1{\hat f}_2$.
Khroshkin, Tolstoy \cite{KT} and Yamane \cite{Y} showed that
\begin{equation}
{\cal R}_q={\hat {\cal R}}_qK_q,
\end{equation}
where
\begin{eqnarray}
&&K_q=q^{-h_1\otimes h_2-h_2\otimes h_1-2h_2\otimes h_2},\nonumber\\
&&{\hat {\cal
R}}_q=\exp_q\left[\left(q-q^{-1}\right)q^{-h_1/2}{\hat e}_1\otimes
{\hat
f}_1q^{h_1/2}\right]\exp_q\left[-\left(q-q^{-1}\right)q^{-h_1/2-h_2/2}e_3\otimes
{\hat f}_3q^{h_1/2+h_2/2}\right]\times
\nonumber\\
&&\phantom{{\cal {\hat
R}}_q=}\exp_q\left[-\left(q-q^{-1}\right)q^{-h_2/2}{\hat
e}_2\otimes {\hat f}_2q^{h_2/2}\right]
\end{eqnarray}
and $\exp_q\left(x\right)=\sum_{n\geq 0}x^n/\left(n\right)_q!$,
$\left(n\right)_q!=\left(1\right)_q\left(2\right)_q\cdots
\left(n\right)_q$, $\left(n\right)_q=\frac{1-q^k}{1-q}$. In the
(fund.)$\otimes$(arbitrary) representation, the $R$-matrix (3)
takes the following form:
\begin{equation}
\left.R_{\sf q}\right|_{\left(fund.\otimes arb.\right)}
=\pmatrix{q^{-h_2}& &A&&B\cr &&&&\cr 0&&q^{-h_1-h_2}&&C\cr &&&&\cr
0&&0&&q^{-h_1-2h_2}},
\end{equation}
where
\begin{eqnarray}
&&A=\left(q-q^{-1}\right)q^{-1/2}{\hat f}_1q^{-h_1/2-h_2},\nonumber\\
&&B=-\left(q-q^{-1}\right)^2q^{-1}{\hat f}_1q^{h_1/2}{\hat
f}_2q^{-h_1-3h_2/2}-\left(q-q^{-1}\right)q^{-1/2}{\hat
f}_3q^{-h_1/2-3h_2/2},
\nonumber\\
&&C=\left(q-q^{-1}\right)q^{-1/2}{\hat f}_2q^{-h_1-3h_2/2},
\end{eqnarray}

2. The $R_{\sf h}$-matrix in the (fund.)$\otimes$(arbitrary)
representation is obtained, from (5), as follows:
\begin{eqnarray}
\left.R_{\sf h}\right|_{\left(fund.\otimes arb.\right)}&=&
\lim_{q\rightarrow 1} \pmatrix{{\sf G}^{-1}& &-\frac{\sf
h}{q-1}{\sf G}^{-1}&&0\cr &&&&\cr 0&&{\sf G}^{-1}&&0\cr &&&&\cr
0&&0&&{\sf G}^{-1}} \left.R_{\sf q}\right|_{\left(fund.\otimes
arb.\right)}\pmatrix{{\sf G}& &\frac{\sf h}{q-1}{\sf G}&&0\cr
&&&&\cr 0&&{\sf G}&&0\cr &&&&\cr
0&&0&&{\sf G}}\nonumber\\
&=& \lim_{q\rightarrow 1} \pmatrix{{\sf G}^{-1}q^{-h_2}{\sf G}&
&\alpha &&\beta \cr &&&&\cr 0&&{\sf G}^{-1}q^{-h_1-h_2}{\sf
G}&&\gamma\cr &&&&\cr 0&&0&&{\sf G}^{-1}q^{-h_1-2h_2}},
\end{eqnarray}
where
\begin{eqnarray}
&&\alpha =\frac{\sf h}{q-1}\left({\sf G}^{-1}q^{-h_2}{\sf G}-{\sf
G}^{-1}q^{-h_1-h_2}{\sf G}\right)+
\left(q-q^{-1}\right)q^{-1/2}{\sf G}^{-1}{\hat f}_1q^{-h_1/2-h_2}{\sf G},\nonumber\\
&&\beta=-\left(q-q^{-1}\right)^2q^{-1}{\sf G}^{-1}{\hat
f}_1q^{h_1/2}{\hat f}_2q^{-h_1-3h_2/2}{\sf G}-
\left(q-q^{-1}\right)q^{-1/2}{\sf G}^{-1}{\hat f}_3q^{-h_1/2-3h_2/2}{\sf G}+\nonumber\\
&&\phantom{\beta=}\frac{\sf
h}{q-1}\left(q-q^{-1}\right)q^{-1/2}{\sf G}^{-1}{\hat f}_2q^{-h_1-3h_2/2}{\sf G}\nonumber\\
&&\gamma=-\left(q-q^{-1}\right)q^{-1/2}{\sf G}^{-1}{\hat
f}_2q^{-h_1-3h_2/2}{\sf G}
\end{eqnarray}
and ${\sf G}=E_q\left(\frac{\sf h}{q-1}{\hat e}_1\right)$,
$E_q\left(x\right)=\sum_{n\geq 0}x^n/[n]!$, $[n]!=[n]\cdots [1]$,
$[n]=\frac{q^n-q^{-n}}{q-q^{-n}}$.

3. Defining
\begin{equation}
t^{(\alpha)}=E_q^{-1}\left(\frac{\sf h}{q-1}{\hat
e}_1\right)E_q\left(q^\alpha\frac{\sf h}{q-1}{\hat
e}_1\right),\qquad t^{(0)}=1,
\end{equation}
we obtain the following properties:
\begin{eqnarray}
&&E_q^{-1}\left(\frac{{\sf h}{\hat e}_1}{q-1}\right)q^{\alpha
h_1/2}E_q\left(\frac{{\sf h}{\hat
e}_1}{q-1}\right)=t^{(\alpha)}q^{\alpha
h_1/2},\\
&&t^{(\alpha+\beta)}q^{(\alpha+\beta)h_1/2}=t^{(\alpha)}q^{\alpha
h_1/2}t^{(\beta)}q^{\beta h_1/2},\\
&&E_q^{-1}\left(\frac{{\sf h}{\hat e}_1}{q-1}\right)q^{\beta
h_2}E_q\left(\frac{{\sf h}{\hat
e}_1}{q-1}\right)=t^{(-\beta)}q^{\beta
h_2},\\
&&E_q^{-1}\left(\frac{{\sf h}{\hat e}_1}{q-1}\right){\hat
f}_1E_q\left(\frac{{\sf h}{\hat e}_1}{q-1}\right)={\hat
f}_1-\frac{\sf
h}{\left(q-1\right)\left(q-q^{-1}\right)}\left(t^{(1)}q^{
h_1}-t^{(-1)}q^{-h_1}\right),\\
&&E_q^{-1}\left(\frac{{\sf
h}{\hat e}_1}{q-1}\right){\hat f}_2E_q\left(\frac{{\sf h}{\hat e}_1}{q-1}\right)={\hat f}_2,\\
&&E_q^{-1}\left(\frac{{\sf h}{\hat e}_1}{q-1}\right){\hat
f}_3E_q\left(\frac{{\sf h}{\hat e}_1}{q-1}\right)={\hat
f}_3+\frac{{\sf h}q}{q-1}t^{(1)}{\hat f}_2q^{ h_1}.
\end{eqnarray}

4. Let us introduce the following generator:
\begin{equation}
T=\lim_{q\rightarrow 1}t^{(1)}.
\end{equation}
From (11) it is evident that
\begin{equation}
\lim_{q\rightarrow 1}t^{(\alpha)}=T^\alpha.
\end{equation}
To obtain a closed form of $T$, we proceed as follows: we left and
right multiply the commutation relation
$q^{h_1}-q^{-h_1}=\left(q-q^{-1}\right)\left({\hat e}_1{\hat
f}_1-{\hat f}_1{\hat e}_1\right)$ by ${\sf G}^{-1}$ and ${\sf G}$
respectively. After simple calculations, we reach to
\begin{equation}
t^{(2)}q^{h_1}-t^{(-2)}q^{-h_1}=q^{h_1}-q^{-h_1}+{\sf
h}\left(q+1\right)\left[t^{(1)}{\hat e}_1q^{h_1}+q^{-h_1}{\hat
e}_1t^{(-1)}\right],
\end{equation}
which yields, when $q\longrightarrow 1$, to
\begin{equation}
T^2-T^{-2}=2{\sf h}\left(T+T^{-1}\right){\hat e}_1\qquad
\Rightarrow\qquad  T-T^{-1}=2{\sf h}{\hat e}_1.
\end{equation}
Finally, we obtain
\begin{eqnarray}
&&T^{\pm 1}=\pm {\sf h}{\hat e}_1+\sqrt{1+{\sf h}^2{\hat e}_1^2}.
\end{eqnarray}

5. We turn now to (7). It is easy to verify that
\begin{eqnarray}
&&\lim_{q\rightarrow 1}{\sf {\sf G}}^{-1}q^{-h_2}{\sf {\sf
G}}=\lim_{q\rightarrow
1}t^{(1)}q^{-h_2}=T,\nonumber\\
&&\lim_{q\rightarrow 1}{\sf G}^{-1}q^{-h_1-h_2}{\sf
G}=\lim_{q\rightarrow
1}t^{(-1)}q^{-h_1-h_2}=T^{-1},\nonumber\\
&&\lim_{q\rightarrow
1}{\sf G}^{-1}q^{-h_1-2h_2}{\sf G}=1, \nonumber\\
&&\lim_{q\rightarrow 1}\gamma=-\lim_{q\rightarrow
1}\left(q-q^{-1}\right)q^{-1/2}{\hat
f}_2t^{(-1/2)}q^{-h_1-3h_2/2}=0,
\nonumber\\
&&\lim_{q\rightarrow 1}\alpha=\lim_{q\rightarrow
1}\left(q-q^{-1}\right)q^{-1/2}{\hat f}_1+\frac{{\sf
h}t^{(1)}}{q-1}\left(q^{-h_2}-q^{-h_1/2-h_2-1/2}\right)
-\frac{{\sf
h}t^{(-1)}}{q-1}\left(q^{-h_1-h_2}-q^{-3h_1/2-h_2-1/2}\right)\nonumber\\
&&\phantom{\lim_{q\rightarrow 1}\alpha}=-\frac{\sf
h}2\left(T+T^{-1}\right)h_1+\frac{\sf
h}2\left(T-T^{-1}\right)\nonumber\\
&& \phantom{\lim_{q\rightarrow 1}\alpha}\equiv -{\sf
h}H_1+\frac{\sf h}2\left(T-T^{-1}\right),\nonumber\\
&&\lim_{q\rightarrow 1}\beta=-\lim_{q\rightarrow
1}\left(q-q^{-1}\right)^2q^{-1}\left[{\hat f}_1-\frac{{\sf
h}}{\left(q-1\right)\left(q-q^{-1}\right)}\left(t^{(1)}q^{h_1}-t^{(-1)}q^{-h_1}\right)\right]
t^{(1)}q^{h_1^/2}{\hat f}_2t^{(-1/2)}q^{-h_1-3h_2/2}\nonumber\\
&&\phantom{\lim_{q\rightarrow
1}\beta=}-\left(q-q^{-1}\right)q^{-1/2}\left[f_3+\frac{{\sf
h}q}{q-1}t^{(1)}q^{h_1}\right]t^{(1/2)}q^{-h_1/2-3h_2/2}\nonumber\\
&&\phantom{\lim_{q\rightarrow 1}\beta=}+\frac{{\sf
h}\left(q-q^{-1}\right)}{q-1}q^{-1/2}{\hat
f}_2t^{(-1/2)}q^{-h_1-3h_2/2}\nonumber\\
&&\phantom{\lim_{q\rightarrow 1}\beta}=2{\sf
h}\left(T-T^{-1}\right)T^{1/2}f_2-2{\sf h}T^{3/2}f_2+2{\sf
h}T^{-1/2}f_2=0.
\end{eqnarray}
Finally, we obtain
\begin{eqnarray}
\left.R_{\sf h}\right|_{\left(fund.\otimes fund.\right)}&=&
\pmatrix{T& &-{\sf h}H_1+\frac{\sf h}2\left(T-T^{-1}\right)&&0\cr
&&&&\cr 0&&T^{-1}&&0\cr &&&&\cr 0&&0&&1}.
\end{eqnarray}

\end{document}